\documentclass{amsart}
\usepackage{amsthm}                 %Definitionen, Theoreme etc.
\usepackage{amsmath}                %siehe userguide, z.B. boxed{}, text-Befehl etc.
\usepackage{amsfonts,amssymb}       %ein paar zusätzliche Symbole
\usepackage{graphicx}               %Bilder
\usepackage{float}                  %Bilder/Positionierung
\usepackage{rotating}               %Bilder/Positionierung
\restylefloat{figure}               %Bilder/Positionierung
\usepackage{mathrsfs}               %mathscr-Buchstaben
\usepackage{fancyhdr}               %Kopfzeile
\usepackage{hyperref}               %Referenzen anklicken
\usepackage[numbers,sort&compress]{natbib}
\usepackage{color}
\theoremstyle{plain}
\newtheorem{lem}{Lemma}

\newtheorem{thm}[lem]{Theorem}

\theoremstyle{definition}

\newcommand{\R}{\mathbb R}

\newcommand{\Z}{\mathbb Z}
\newcommand{\N}{\mathbb N}

\renewcommand{\L}{\mathcal L}

\renewcommand{\d}{\,\text{\rm d}}
\newcommand{\dw}{\text{\rm d}}

\renewcommand{\S}{\mathbb S}
\renewcommand{\phi}{\varphi}

\newcommand{\norm}[1]{\left|\!\left|#1\right|\!\right|}

\newcommand{\eps}{\varepsilon}

\newcommand{\ska}[2]{\left\langle #1,#2\right\rangle}

\newcommand{\set}[2]{\left\{#1;\;#2\right\}}

\newcommand{\bea}{\begin{eqnarray}}
\newcommand{\eea}{\end{eqnarray}}
\newcommand{\beq}{\begin{equation}}
\newcommand{\eeq}{\end{equation}}
\renewcommand{\phi}{\varphi}

\renewcommand{\autoref}[1]{\text{Eq.}~\eqref{#1}}

\newcommand{\A}{\mathcal A}
\newcommand{\B}{\mathcal B}
\newcommand{\RR}{\mathcal R}
\newcommand{\M}{\mathcal M}
\newcommand{\T}{\mathcal T}
\newcommand{\calS}{\mathcal S}
\begin{document}
\title{Necrotic tumor growth: an analytic approach}
\author{Martin Kohlmann}
\address{Peter L. Reichertz Institute for Medical Informatics, University of Braunschweig, D-38106 Braunschweig, Germany}
\email{martin.kohlmann@plri.de}
\begin{abstract}
The present paper deals with a free boundary problem modeling the growth process of necrotic multi-layer tumors.
We prove the existence of flat stationary solutions and determine the linearization of our model at such an equilibrium.
Finally, we compute the solutions of the stationary linearized problem and comment on bifurcation.
\end{abstract}
\maketitle
\tableofcontents
\section{Introduction}\label{sec_intro}
Mathematical models for tumor growth have been considered with regularity the the applied sciences literature in recent years.
From the mathematical point of view, free boundary models are of particular interest: In these models, a tumor cell at time $t\geq 0$ is identified with
an open domain $\Omega(t)\subset\R^n$, for some $n\geq 1$, with initial configuration $\Omega(0)=\Omega_0$. For simplicity, it is assumed in many models that the growth process of the tumor is controlled by only two quantities: the concentration of nutrient (e.g., glucose or oxygen), denoted as $\sigma(x,t)$, and an internal pressure $p(x,t)$, which both have to solve an elliptic problem on the time-dependent and unknown domain $\Omega(t)$, with suitable conditions on the free boundary $\partial\Omega(t)$. Finally, an evolution equation for the free boundary $\partial\Omega(t)$ is needed, and usually it is derived from a simple application of Darcy's law pertaining to the fact that the tumor behaves as an incompressible ideal fluid.

In many publications dealing with free boundary problems for tumor growth, the domain $\Omega(t)$ is assumed to be spherically symmetric and $n=1$, cf. the seminal papers \cite{BC95,FR99,WK97}. The present work is innovative for the following three reasons:
\begin{itemize}
\item We are looking at strip-shaped tumors: Lately, biologists have discovered that a peculiar kind of in vitro tumors can by cultivated by a special tissue culture technique, cf.~\cite{KSoH04,KCM99,M97}. In our model, we imagine the tumor to start growing from the flat bottom of a Petri dish.
\item An additional feature of our model is that we distinguish between a necrotic core, localized at the bottom of the Petri dish, and a non-necrotic shell which is lying above. In consequence, our problem has two free boundaries confining a time-dependent domain on which we study elliptic problems for nutrient and pressure.
\item Finally, we present a two-dimensional model (i.e., $n=2$).
\end{itemize}
We refer the reader to \cite{CuiE09,ZECui08b}, where the authors explain a sophisticated approach to the growth of non-necrotic multi-layer tumors, and \cite{EMM10}, where spherically symmetric necrotic tumor cells are studied. Based on the model assumptions in \cite{CuiE09,EMM10}, we now present the following problem:

Let $\S=\R/2\pi\Z$ and consider two positive time-dependent functions $\rho_1<\rho_2$ on $\S$. Let furthermore
$$\Omega_{\rho_1,\rho_2}(t)=\set{(x,y)\in\R^2}{x\in\S,\,\rho_1(t,x)<y<\rho_2(t,x)}$$
with the boundary components
$$
\Gamma_{\rho_i}(t)=\set{(x,y)\in\R^2}{y=\rho_i(t,x)},\quad i=1,2.
$$
The outward unit normal of $\Gamma_{\rho_i}(t)$ with respect to $\Omega_{\rho_1,\rho_2}(t)$ is denoted by $\nu_i$, for $i=1,2$. We obtain $\nu_1$ and $\nu_2$ by computing the gradients of the functions $N_i(x,y)=y-\rho_i(x)$:
$$\nu_1=-\frac{\nabla N_1}{|\nabla N_1|}\quad\text{and}\quad\nu_2=\frac{\nabla N_2}{|\nabla N_2|}.$$
We will write $n_1=-\nabla N_1$ and $n_2=\nabla N_2$. Let $\kappa_i$ denote the
curvature of $\Gamma_{\rho_i}(t)$. It is well known that $\kappa_i$ can be computed explicitly using the formula
$$\kappa_i=-\frac{\rho_{ixx}}{(1+\rho_{ix}^2)^{3/2}},\quad i=1,2.$$
Let $0<\rho_{1,0}<\rho_{2,0}$ be periodic functions on $\R$ so that $\rho_i(x,0)=\rho_{i,0}(x)$. The nutrient $\sigma$ should satisfy a stationary diffusion equation. Furthermore, we assume that there is a constant supply $\bar\sigma>0$ of nutrient on $\Gamma_{\rho_2}(t)$ and that the normal derivative of $\sigma$ vanishes on $\Gamma_{\rho_1}(t)$. Next, the Laplacian of the pressure is proportional to the difference $\sigma-\tilde\sigma$, with proportionality factor $-\mu$; here $\tilde\sigma$ and $\mu$ are positive parameters. The reason for this assumption is that if $\sigma<\tilde\sigma$, then the tumor volume locally decreases, whereas the tumor grows in regions where $\sigma>\tilde\sigma$. The boundary conditions for the pressure are the so-called Laplace-Young conditions: We assume that the pressure on $\Gamma_{\rho_{2}}(t)$ is proportional to $\kappa_2$; the proportionality constant is the surface tension coefficient $\gamma_2>0$. Similarly, we have $p=\gamma_1\kappa_1-c$, with $\gamma_1>0$ and $c$ a positive constant.
We also assume that
\beq\bar\sigma>\tilde\sigma\label{bartilde}\eeq
to have a reasonable long-time behavior, as explained in \cite{CuiE09}.
Finally, the normal velocity of the boundary components is equal to the cell movement velocity in the direction $n_2$ on $\Gamma_{\rho_2}$ and $-n_1$ on $\Gamma_{\rho_1}$ respectively. This yields two evolution equations for the moving boundaries.
\begin{figure}[H]
\begin{center}
\includegraphics[width=8cm]{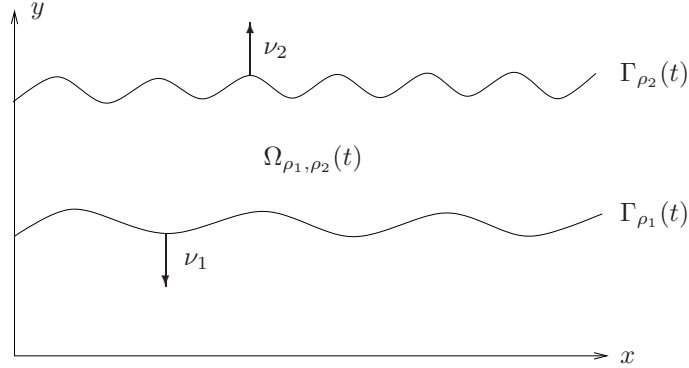}
\end{center}
\caption{A free boundary problem for the growth of multi-layer tumors with a necrotic core.}
\end{figure}
\begin{picture}(0,0)
\put(120,95){$\nu_1$}
\put(150,175){$\nu_2$}
\put(113,107){\vector(0,-1){20}}
\put(145,167){\vector(0,1){20}}
\put(150,133){$\Omega_{\rho_1,\rho_2}(t)$}
\put(285,112){$\Gamma_{\rho_1}(t)$}
\put(285,165){$\Gamma_{\rho_2}(t)$}
\put(62,190){$y$}
\put(285,58){$x$}
\end{picture}
Our mathematical model is given by the following system of equations:
\beq
\left\{
\begin{array}{rclcl}
  \Delta\sigma                         & = & \sigma                             & \text{in}  & \Omega_{\rho_1,\rho_2}(t),\\
  \Delta p                             & = & -\mu(\sigma-\tilde\sigma)          & \text{in}  & \Omega_{\rho_1,\rho_2}(t), \\
  \frac{\partial\sigma}{\partial n_1}  & = & 0                                  & \text{on}  & \Gamma_{\rho_1}(t),\\
  p                                    & = & \gamma_1\kappa_1 - c               & \text{on}  & \Gamma_{\rho_1}(t),\\
  \sigma                               & = & \bar\sigma                         & \text{on}  & \Gamma_{\rho_2}(t),\\
  p                                    & = & \gamma_2\kappa_2                   & \text{on}  & \Gamma_{\rho_2}(t),\\
  \rho_{1,t}                           & = &  \frac{\partial p}{\partial n_1}   & \text{on}  & \Gamma_{\rho_1}(t),\\
  \rho_{2,t}                           & = & -\frac{\partial p}{\partial n_2}   & \text{on}  & \Gamma_{\rho_2}(t),\\
  \rho_1                               & = & \rho_{1,0}                         & \text{for} & t=0,\\
  \rho_2                               & = & \rho_{2,0}                         & \text{for} & t=0.\\
\end{array}
\right.
\label{problem}
\eeq
A solution to \eqref{problem} is a tupel $(\sigma(x,y,t),p(x,y,t),\rho_{1}(x,t),\rho_2(x,t))$, where $t\in[0,T]$, $T>0$, $x\in\S$ and $\rho_1(x,t)\leq y\leq\rho_2(x,t)$, so that all equations of \eqref{problem} are satisfied pointwise.
In the following sections, we will only discuss the stationary version of \eqref{problem}. It is given by the following system of equations
\beq
\left\{
\begin{array}{rclcl}
  \Delta\sigma                         & = & \sigma                    & \text{in}  & \Omega_{\rho_1,\rho_2},\\
  \Delta p                             & = & -\mu(\sigma-\tilde\sigma) & \text{in}  & \Omega_{\rho_1,\rho_2}, \\
  \frac{\partial\sigma}{\partial n_1}  & = & 0                         & \text{on}  & \Gamma_{\rho_1},\\
  \frac{\partial p}{\partial n_1}      & = & 0                         & \text{on}  & \Gamma_{\rho_1},\\
  p                                    & = & \gamma_1\kappa_1 - c      & \text{on}  & \Gamma_{\rho_1},\\
  \sigma                               & = & \bar\sigma                & \text{on}  & \Gamma_{\rho_2},\\
  \frac{\partial p}{\partial n_2}      & = & 0                         & \text{on}  & \Gamma_{\rho_2},\\
  p                                    & = & \gamma_2\kappa_2          & \text{on}  & \Gamma_{\rho_2},\\
\end{array}
\right.
\label{stationary}
\eeq
where $(\sigma,p,\rho_1,\rho_2)=(\sigma(x,y),p(x,y),\rho_1(x),\rho_2(x))$.

The layout of this paper is as follows: We first prove the existence of flat stationary solutions, i.e., solutions
$(\sigma,p,\rho_1,\rho_2)$ with $\sigma=\sigma(y)$, $p=p(y)$ and with constants $0<\rho_1<\rho_2$. Precisely, for any given set of
positive values $\bar\sigma>\tilde\sigma$ and $\mu$ there is a
flat stationary solution which is unique up to a shift in the $y$-direction. Next, we linearize the system \eqref{stationary} at such an equilibrium and use Fourier expansions to obtain solutions of the linearized system. The calculations carried out in the following sections especially show how to choose $c$ and the surface tension coefficients $\gamma_1$ and $\gamma_2$ to obtain non-trivial solutions. In an outlook, we present the bifurcation problem associated with the model~\eqref{stationary}.
\\[.1cm]

\emph{Acknowledgement.} The author thanks the anonymous referee for asking about the bifurcation problem associated with the model of the paper at hand which led to an additional chapter compared to the initially submitted version.
\section{Flat stationary solutions}
Let $(\sigma_*,p_*,\rho_{1*},\rho_{2*})$ be a flat stationary solution of the problem \eqref{stationary}, i.e., we have that
\beq
\left\{
\begin{array}{rcl}
  \sigma_*''             & = & \sigma_*,                     \\
  p_*''                  & = & -\mu(\sigma_*-\tilde\sigma),  \\
  \sigma_*(\rho_{2*})    & = & \bar\sigma,                   \\
  \sigma_*'(\rho_{1*})   & = & 0,                            \\
  p_*(\rho_{1*})         & = & -c,                           \\
  p_*(\rho_{2*})         & = & 0,                            \\
  p_*'(\rho_{1*})        & = & 0,                            \\
  p_*'(\rho_{2*})        & = & 0.                            \\
\end{array}
\right.
\nonumber
\eeq
We first solve the subproblem for the nutrient concentration and find that
\beq\label{sigmastar}\sigma_*(y)=\bar\sigma\frac{\cosh y-\tanh\rho_{1*}\sinh y}{\cosh\rho_{2*}-\tanh\rho_{1*}\sinh\rho_{2*}}\eeq
is the unique solution of
\beq
\left\{
\begin{array}{rcl}
  \sigma_*''             & = & \sigma_*,                     \\
  \sigma_*(\rho_{2*})    & = & \bar\sigma,                   \\
  \sigma_*'(\rho_{1*})   & = & 0.                            \\
\end{array}
\right.
\nonumber
\eeq
Next we consider the boundary value problem
\beq
\left\{
\begin{array}{rcl}
  p_*''                  & = & -\mu(\sigma_*-\tilde\sigma),  \\
  p_*(\rho_{2*})         & = & 0,                            \\
  p_*'(\rho_{1*})        & = & 0,                            \\
\end{array}
\right.
\nonumber
\eeq
which has the unique solution
\beq\label{pstar}p_*(y)=\mu(\bar\sigma-\sigma_*)+\mu\tilde\sigma\rho_{1*}(\rho_{2*}-y)+\frac{1}{2}\mu\tilde\sigma(y^2-\rho_{2*}^2).\eeq
Since we must demand $p_*'(\rho_{2*})=0$, we get the condition
\beq\frac{\tilde\sigma}{\bar\sigma}=\frac{\tanh(\rho_{2*}-\rho_{1*})}{\rho_{2*}-\rho_{1*}}.\label{condition1}\eeq
Letting $\rho_{2*}-\rho_{1*}=\delta$, it follows from $0<\tilde\sigma<\bar\sigma$ that \autoref{condition1} has a unique solution $\delta\in(0,\infty)$.
Next the constraint $p_*(\rho_{1*})=-c$ results in the condition
\beq\mu\bar\sigma-\frac{1}{2}\mu\tilde\sigma\delta^2=\mu\bar\sigma
\frac{\cosh\rho_{1*}-\tanh\rho_{1*}\sinh\rho_{1*}}{\cosh(\rho_{1*}+\delta)-\tanh\rho_{1*}\sinh(\rho_{1*}+\delta)}-c.\label{condition2}\eeq
In view of the addition theorems for hyperbolic functions and the relation \eqref{condition1}, \autoref{condition2} can be simplified to
\beq\label{condc}c=\frac{\mu\bar\sigma}{\cosh\delta}\left(1-\cosh\delta+\frac{1}{2}\delta\sinh\delta\right).\eeq
Since the term in brackets is positive for any $\delta>0$, we have $c>0$. Moreover, there is no condition on $\rho_{1*}$, so that we obtain a flat stationary solution of \eqref{problem} for any fixed $\rho_{1*}>0$. This provides a proof of the following theorem.
\thm\label{thm1} Fix $\mu,\bar\sigma,\tilde\sigma,\gamma_1,\gamma_2>0$ and assume that \eqref{bartilde} holds true. Let $\delta$ be the solution of $\frac{\tilde\sigma}{\bar\sigma}=\frac{\tanh\delta}{\delta}$ and define $c$ according to \autoref{condc}. Then there is a one-parameter family of flat stationary solutions $(\sigma_*,p_*,\rho_{1*},\rho_{1*}+\delta)$ to \autoref{problem}, where $\sigma_*$ and $p_*$ are given by \autoref{sigmastar} and \autoref{pstar}, respectively.
\endthm\rm
\section{The linearized problem and its solutions}\label{sec_3}
A standard technique to tackle moving boundary problems is to transform the problem under consideration to a problem on a fixed (and preferably simple) reference domain, to solve the problem on the reference domain and to transform its solutions back to obtain solutions of the original problem.

Assume that $\rho_{1},\rho_{2}\in C^2_+(\S)$ and let $\vartheta(x,y)=\frac{y-\rho_1(x)}{\rho_2(x)-\rho_1(x)}$, so that the map $\psi\colon(x,y)\mapsto (x',y'):=(x,\vartheta(x,y))$ establishes a $C^2$-diffeomorphism $\Omega_{\rho_1,\rho_2}\to\S\times(0,1)$. The strip $\Omega=\S\times(0,1)$ will be our reference domain, with the boundary components $\Gamma_1=\S\times\{0\}$ and $\Gamma_2=\S\times\{1\}$. In particular, we have that $\Gamma_i\simeq\S$ and that $\psi(\Gamma_{\rho_i})=\Gamma_i$, for $i=1,2$.
\begin{figure}[H]
\begin{center}
\includegraphics[width=12cm]{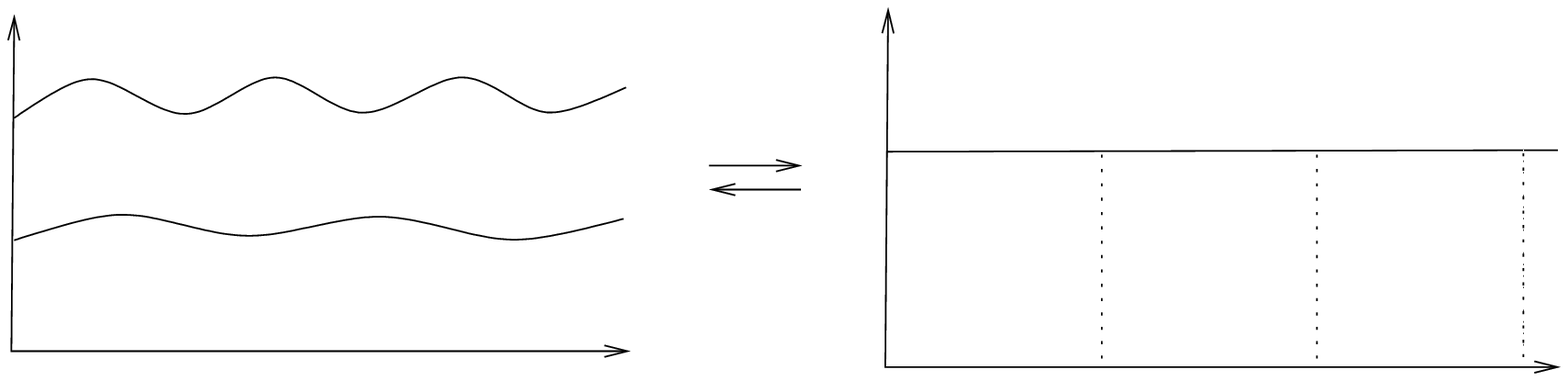}
\end{center}
\caption{Transformation of the tumor domain onto a fixed reference domain.}
\end{figure}
\begin{picture}(0,0)
\put(155,82){$\psi^{-1}$}
\put(157,104){$\psi$}
\put(80,76){$\Gamma_{\rho_1}$}
\put(80,122){$\Gamma_{\rho_2}$}
\put(55,96){$\Omega_{\rho_1,\rho_2}$}
\put(140,54){$x$}
\put(05,130){$y$}
\put(344,51){$x'$}
\put(197,130){$y'$}
\put(257,74){$\Omega$}
\put(300,105){$\Gamma_2$}
\put(300,58){$\Gamma_1$}
\end{picture}
Let us further introduce the operators
$$\mathcal A(\rho_1,\rho_2)u=[\Delta (u\circ\psi)]\circ\psi^{-1}\quad\text{and}\quad\mathcal B_i(\rho_1,\rho_2)u=\ska{\text{tr}_i[\nabla (u\circ\psi)]}{n_i}\circ\psi^{-1},$$
where $\text{tr}_i$ denotes the trace with respect to $\Gamma_{\rho_i}$, for $i=1,2$, and $u\in C^2(\overline\Omega)$. A straightforward computation shows that
\bea\mathcal A(\rho_1,\rho_2)u&=&u_{x'x'}-2u_{x'y'}\frac{y'(\rho_2'-\rho_1')+\rho_1'}{\rho_2-\rho_1}
+u_{y'y'}\frac{1+[y'(\rho_2'-\rho_1')+\rho_1']^2}{(\rho_2-\rho_1)^2}\nonumber\\
&&\quad+u_{y'}\left(2\frac{\rho_2'-\rho_1'}{(\rho_2-\rho_1)^2}[y'(\rho_2'-\rho_1')+\rho_1']-\frac{y'(\rho_2''-\rho_1'')+\rho_1''}{\rho_2-\rho_1}\right)\nonumber
\eea
and
\bea
\mathcal B_1(\rho_1,\rho_2)u &=& \rho_1'u_{x'}|_{y'=0}-\frac{\rho_1'^2+1}{\rho_2-\rho_1}u_{y'}|_{y'=0},\nonumber\\
\mathcal B_2(\rho_1,\rho_2)u &=& -\rho_2'u_{x'}|_{y'=1}+\frac{\rho_2'^2+1}{\rho_2-\rho_1}u_{y'}|_{y'=1}.\nonumber
\eea
Hence the transformed problem reads
\beq\label{trsf}
\left\{
\begin{array}{rclcl}
  \A(\rho_1,\rho_2)\sigma              & = & \sigma                             & \text{in}  & \Omega\times[0,T],  \\
  \A(\rho_1,\rho_2)p                   & = & -\mu(\sigma-\tilde\sigma)          & \text{in}  & \Omega\times[0,T],  \\
  \B_1(\rho_1,\rho_2)\sigma            & = & 0                                  & \text{on}  & \Gamma_1\times[0,T],\\
  p                                    & = & \gamma_1\kappa_1 - c               & \text{on}  & \Gamma_1\times[0,T],\\
  \sigma                               & = & \bar\sigma                         & \text{on}  & \Gamma_2\times[0,T],\\
  p                                    & = & \gamma_2\kappa_2                   & \text{on}  & \Gamma_2\times[0,T],\\
  \rho_{1,t}                           & = &  \B_1(\rho_1,\rho_2)p              & \text{on}  & \Gamma_1\times[0,T],\\
  \rho_{2,t}                           & = & -\B_2(\rho_1,\rho_2)p              & \text{on}  & \Gamma_2\times[0,T],\\
  \rho_1                               & = & \rho_{1,0}                         & \text{for} & t=0,\\
  \rho_2                               & = & \rho_{2,0}                         & \text{for} & t=0.\\
\end{array}
\right.
\eeq
We now pick a flat stationary solution $(\sigma_*,p_*,\rho_{1*},\rho_{2*})$ as obtained in Theorem~\ref{thm1} and, for $\eps>0$, we let
\beq
\left(
  \begin{array}{c}
    \sigma(x',y',t) \\
    p(x',y',t) \\
    \rho_1(x',t) \\
    \rho_2(x',t)
  \end{array}
\right)
=
\left(
  \begin{array}{c}
    \sigma_*(y'\delta+\rho_{1*}) \\
    p_*(y'\delta+\rho_{1*}) \\
    \rho_{1*} \\
    \rho_{2*}
  \end{array}
\right)
+\eps
\left(
  \begin{array}{c}
    \Sigma(x',y',t) \\
    P(x',y',t) \\
    r(x',t) \\
    s(x',t)
  \end{array}
\right),
\label{eps}\eeq
where $\delta=\rho_{2*}-\rho_{1*}>0$. Our regularity assumption on the new unknowns $(\Sigma,P,r,s)$ is that $\Sigma,P\in C([0,T];\,C^2(\overline\Omega))$ and $r,s\in C^1([0,T];\,C^2(\S))$. We also introduce the second order linear differential operator $b_{r,s}$ on $C^2([0,1])$ given by
$$b_{r,s}(u)=\frac{2}{\delta}(s-r)u''(y'\delta+\rho_{1*})+[y'(s''-r'')+r'']u'(y'\delta+\rho_{1*}).$$
The linearization of the problem \eqref{trsf} at the flat stationary solution $(\sigma_*,p_*,\rho_{1*},\rho_{2*})$ is obtained by inserting \eqref{eps} into the problem \eqref{trsf} and by differentiating each equation with respect to $\eps$ at $\eps=0$. This yields
\bea
\left\{
\begin{array}{l}
  \Sigma_{x'x'}+\frac{1}{\delta^2}\Sigma_{y'y'}=b_{r,s}(\sigma_*)+\Sigma, \quad\text{in }\Omega\times[0,T],\\
  \Sigma_{y'}(x',0,t)=\Sigma(x',1,t)=0, \quad\text{on }\S\times[0,T],\\[.25cm]
  P_{x'x'}+\frac{1}{\delta^2}P_{y'y'}=b_{r,s}(p_*)-\mu\Sigma, \quad\text{in }\Omega\times[0,T],\\
  P(x',0,t)=-\gamma_1 r''(x',t), \quad\text{on }\S\times[0,T],\\
  P(x',1,t)=-\gamma_2 s''(x',t), \quad\text{on }\S\times[0,T],\\[.25cm]
  r_t(x',t)=-\frac{1}{\delta}P_{y'}(x',0,t), \quad\text{on }\S\times[0,T],\\
  s_t(x',t)=-\frac{1}{\delta}P_{y'}(x',1,t), \quad\text{on }\S\times[0,T],\\
  r(x',0)=r_0(x'), \quad\text{on }\S,\\
  s(x',0)=s_0(x'), \quad\text{on }\S.
\end{array}
\right.
\label{lin}\eea
The stationary version of \eqref{lin} is
\bea
\left\{
\begin{array}{l}
  \Sigma_{x'x'}+\frac{1}{\delta^2}\Sigma_{y'y'}=b_{r,s}(\sigma_*)+\Sigma, \quad\text{in }\Omega, \\
  \Sigma_{y'}(x',0)=\Sigma(x',1)=0, \quad\text{on }\S,\\[.25cm]
  P_{x'x'}+\frac{1}{\delta^2}P_{y'y'}=b_{r,s}(p_*)-\mu\Sigma, \quad\text{in }\Omega,\\
  P(x',0)=-\gamma_1 r'', \quad\text{on }\S,\\
  P(x',1)=-\gamma_2 s'', \quad\text{on }\S,\\[.25cm]
  P_{y'}(x',0)=0, \quad\text{on }\S,\\
  P_{y'}(x',1)=0, \quad\text{on }\S.
\end{array}
\right.
\label{linstat}\eea
To obtain solutions of \eqref{linstat}, we expand $(\Sigma,P,r,s)$ as Fourier series and denote the coefficients with respect to the basis functions $\{1,\cos(kx'),\sin(kx');\;k\in\N\}$ by
\bea
&\{A_0(y'),A_k(y'),B_k(y');\;k\in\N\} & \text{for the variable } \Sigma,\nonumber\\
&\{M_0(y'),M_k(y'),N_k(y');\;k\in\N\} & \text{for the variable } P,\nonumber\\
&\{a_0,a_k,b_k;\;k\in\N\} & \text{for the variable } r,\nonumber\\
&\{c_0,c_k,d_k;\;k\in\N\} & \text{for the variable } s.\nonumber
\eea
We now proceed in the following steps: First, we solve the boundary value problems
\bea
\left\{
\begin{array}{rcl}
  -k^2A_k(y')+\frac{1}{\delta^2}A_k''(y') & = & A_k(y')+f_k(y'), \\
  A_k'(0) & = & 0, \\
  A_k(1)  & = & 0,
\end{array}
\right.
\label{problemAk}
\eea
for $k=0,1,2,\ldots$, and
\bea
\left\{
\begin{array}{rcl}
  -k^2B_k(y')+\frac{1}{\delta^2}B_k''(y') & = & B_k(y')+g_k(y'), \\
  B_k'(0) & = & 0, \\
  B_k(1)  & = & 0,
\end{array}
\right.
\label{problemBk}
\eea
for $k=1,2,\ldots$; here
\begin{align}
f_k(y')&=\frac{2}{\delta}(c_k-a_k)\sigma_*''(y'\delta+\rho_{1*})-k^2[y'(c_k-a_k)+a_k]\sigma_*'(y'\delta+\rho_{1*})\nonumber\\
&=\frac{2\bar\sigma(c_k-a_k)}{\delta\cosh\delta}\cosh(y'\delta)-\frac{k^2\bar\sigma a_k}{\cosh\delta}\sinh(y'\delta)
-\frac{k^2\bar\sigma(c_k-a_k)}{\cosh\delta}y'\sinh(y'\delta),\nonumber\end{align}
where we have used once again the addition theorems for hyperbolic functions.
The $g_k(y')$ are similar; we simply have to replace the $a_k$ with the $b_k$ and the $c_k$ with the $d_k$. It is straightforward to obtain that
$$
A_0(y')=\frac{\bar\sigma(c_0-a_0)}{\cosh\delta}[y'\sinh(y'\delta)-\tanh\delta\cosh(y'\delta)]
$$
and
\begin{align}
A_k(y')&=\left(\frac{\bar\sigma a_k\tanh(\delta\sqrt{1+k^2})}{\sqrt{1+k^2}\cosh\delta}-\frac{\bar\sigma c_k\tanh\delta}{\cosh(\delta\sqrt{1+k^2})}\right)\cosh(\delta\sqrt{1+k^2}y')\nonumber\\
&\quad-\frac{\bar\sigma a_k}{\sqrt{1+k^2}\cosh\delta}\sinh(\delta\sqrt{1+k^2}y')+\frac{\bar\sigma a_k}{\cosh\delta}\sinh(y'\delta)\nonumber\\
&\quad+\frac{\bar\sigma(c_k-a_k)}{\cosh\delta}y'\sinh(y'\delta).\label{solAk}
\end{align}
The $B_k$ emerge from the $A_k$ by exchanging the $a_k$ with the $b_k$ and the $c_k$ with the $d_k$. We next turn our attention to $M_0$ which has to satisfy the following conditions:
\bea
\left\{
\begin{array}{rcl}
  \frac{1}{\delta^2}M_0''(y') & = & -\mu A_0(y')+\frac{2\mu}{\delta}(c_0-a_0)(\tilde\sigma-\bar\sigma\frac{\cosh(y'\delta)}{\cosh\delta}), \\
  M_0'(0)  & = & 0, \\
  M_0'(1)  & = & 0, \\
  M_0(0)   & = & 0, \\
  M_0(1)   & = & 0. \\
\end{array}
\right.
\label{problemM0}
\eea
Explicit calculations show that there is a solution $M_0$ if and only if $c_0=a_0$ and precisely $M_0\equiv 0$.
Next, we plug the solutions $A_k$ and $B_k$ into the problems
\bea
\left\{
\begin{array}{rcl}
  -k^2M_k(y')+\frac{1}{\delta^2}M_k''(y') & = & -\mu A_k(y')+\tilde f_k(y'), \\
  M_k'(0)  & = & 0, \\
  M_k'(1)  & = & 0,
\end{array}
\right.
\label{problemMk}
\eea
and
\bea
\left\{
\begin{array}{rcl}
  -k^2N_k(y')+\frac{1}{\delta^2}N_k''(y') & = & -\mu B_k(y')+\tilde g_k(y'), \\
  N_k'(0)  & = & 0, \\
  N_k'(1)  & = & 0,
\end{array}
\right.
\label{problemNk}
\eea
for $k=1,2,\ldots$; here
\bea\tilde f_k(y')&=&\frac{2\mu}{\delta}(c_k-a_k)\left(\tilde\sigma-\sigma_*(y'\delta+\rho_{1*})\right)\nonumber\\
&&\quad-\mu k^2(y'(c_k-a_k)+a_k)\left(\tilde\sigma\delta y'-\sigma_*'(y'\delta+\rho_{1*})\right)\nonumber\\
&=&\frac{2}{\delta}\mu\tilde\sigma(c_k-a_k)-\mu k^2\tilde\sigma\delta a_k y'-\mu k^2\tilde\sigma\delta(c_k-a_k)y'^2+\frac{\mu k^2\bar\sigma a_k}{\cosh\delta}\sinh(y'\delta)\nonumber\\
&&\quad-\frac{2\mu\bar\sigma(c_k-a_k)}{\delta\cosh\delta}\cosh(y'\delta)+\frac{\mu k^2\bar\sigma(c_k-a_k)}{\cosh\delta}y'\sinh(y'\delta)\nonumber\eea
and $\tilde g_k(y')$ is obtained as before. Again, it is straightforward to derive the solutions
\bea M_k(y') \!&=&\! -\frac{\mu\cosh(y'\delta k)}{\delta k\sinh(\delta k)}\bigg[\frac{\bar\sigma\delta a_k}{\cosh\delta\cosh(\delta\sqrt{1+k^2})}-\tilde\sigma\delta a_k\cosh(\delta k)\nonumber\\
&&\quad+\bar\sigma c_k\delta\sqrt{1+k^2}\tanh\delta\tanh(\delta\sqrt{1+k^2})+\bar\sigma c_k\tanh\delta-\delta\bar\sigma c_k\bigg]\nonumber\\
&&\quad+\frac{\mu\bar\sigma a_k}{\sqrt{1+k^2}\cosh\delta}\sinh(y'\delta\sqrt{1+k^2})\nonumber\\
&&\quad+\mu\bar\sigma\left(\frac{c_k\tanh\delta}{\cosh(\delta\sqrt{1+k^2})}
-\frac{a_k\tanh(\delta\sqrt{1+k^2})}{\sqrt{1+k^2}\cosh\delta}\right)\cosh(y'\delta\sqrt{1+k^2})\nonumber\\
&&\quad-\frac{\mu\tilde\sigma a_k}{k}\sinh(y'\delta k)-\frac{\mu\bar\sigma a_k}{\cosh\delta}\sinh(\delta y')-\frac{\mu\bar\sigma(c_k-a_k)}{\cosh\delta}y'\sinh(\delta y')\nonumber\\
&&\quad+\mu\tilde\sigma\delta a_ky'+\mu\tilde\sigma\delta(c_k-a_k)y'^2,\label{solMk}\eea
and $N_k$ accordingly. Since we are looking for non-trivial solutions, we assume that $(a_k,b_k),(c_k,d_k)\neq (0,0)$ and that $a_kd_k=b_kc_k$. Without loss of generality, we will henceforth suppose that $a_k\neq 0$ and $c_k\neq 0$. Next we may choose the parameters $\gamma_{1}$ and $\gamma_{2}$ so that
$$\gamma_1=\frac{1}{k^2a_k}M_k(0)\quad\text{and}\quad\gamma_2=\frac{1}{k^2c_k}M_k(1);$$
precisely,
\bea\gamma_1\!&=&\!\frac{\mu\bar\sigma}{\delta k^3\sinh(\delta k)}\bigg[\cosh(\delta k)\tanh\delta-\frac{\delta k\tanh(\delta\sqrt{1+k^2})\sinh(\delta k)}{\sqrt{1+k^2}\cosh\delta}\nonumber\\
&&\quad+\frac{c_k}{a_k}\delta\left(1-\frac{\tilde\sigma}{\bar\sigma}+\tanh\delta\frac{k\sinh(\delta k)-\sqrt{1+k^2}\sinh(\delta\sqrt{1+k^2})}{\cosh(\delta\sqrt{1+k^2})}\right)\nonumber\\
&&\quad-\frac{\delta}{\cosh\delta\cosh(\delta\sqrt{1+k^2})}\bigg]\label{gamma1}\eea
and
\bea\gamma_2\!&=&\!\frac{\mu\bar\sigma}{\delta k^3\tanh(\delta k)}\bigg[\frac{a_k}{c_k}\left(\frac{\tanh\delta}{\cosh(\delta k)}-\frac{\delta}{\cosh\delta\cosh(\delta\sqrt{1+k^2})}\right)\nonumber\\
&&\quad+\delta\tanh\delta\left(k\tanh(\delta k)-\sqrt{1+k^2}\tanh(\delta\sqrt{1+k^2})\right)+\delta-\tanh\delta\bigg].
\label{gamma2}\eea
For simplicity, we will now furthermore assume that $a_k=c_k$ for any $k\in\N$. It is easy to see that the terms in brackets standing one underneath the other in the formula for $\gamma_1$ converge to $+\infty$, $\delta-(1+\frac{\delta^2}{2})\tanh\delta$ and zero, respectively, for $k\to\infty$. Similarly, one deduces that the terms in brackets for $\gamma_2$ tend to $\delta-\tanh\delta$ as $k$ approaches $\infty$. It follows that $\gamma_1,\gamma_2>0$ for $k$ sufficient large. Moreover we have that $\gamma_1,\gamma_2\to 0$, as $k\to\infty$. The figure below shows the functions $\gamma_1(k)$ and $\gamma_2(k)$ for a particular choice of the parameters $\bar\sigma,\tilde\sigma$ and $\mu$.
\begin{figure}[H]
\begin{center}
\includegraphics[width=8cm]{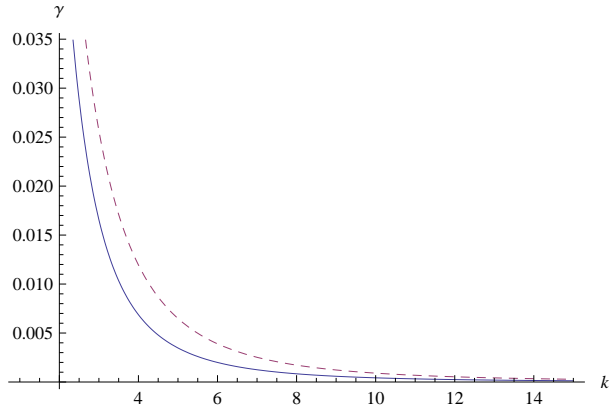}
\end{center}
\caption{The coefficients $\gamma_{1}$ and $\gamma_{2}$ (dashed) for $\tilde\sigma=\frac{1}{2}\bar\sigma=\mu=1$.}
\end{figure}
Note that we have $\sum_k|a_k|^2<\infty$ and $\sum_k|c_k|^2<\infty$, so that $a_k$ and $c_k$ are null sequences, but in general, it is difficult to obtain control of the ratios $\frac{c_k}{a_k}$ and $\frac{a_k}{c_k}$, which is why we are working with the additional assumption $c_k=a_k$ here.

Our result can be formulated as follows.
\thm\label{thm_main} Pick $k\in\N$ and let $a_k,b_k,c_k,d_k$ be real constants satisfying the relations
$$a_k^2+b_k^2\neq 0,\quad c_k^2+d_k^2\neq 0,\quad a_kd_k=b_kc_k.$$
The linearized problem \eqref{linstat} has a nontrivial solution $(\Sigma,P,r,s)$ if and only if $\gamma_1$ and $\gamma_2$ are as in \eqref{gamma1} and \eqref{gamma2} and this solution is given by
\begin{align}
\Sigma(x',y')&=A_k(y')\cos(kx')+B_k(y')\sin(kx'),\nonumber\\
P(x',y')     &=M_k(y')\cos(kx')+N_k(y')\sin(kx'),\nonumber\\
r(x')        &=a_k\cos(kx')+b_k\sin(kx'),\nonumber\\
s(x')        &=c_k\cos(kx')+d_k\sin(kx');\nonumber
\end{align}
the coefficients $A_k$ and $M_k$ are as in \eqref{solAk} and \eqref{solMk} and $B_k$ and $N_k$ are obtained by exchanging $a_k$ with $b_k$ and $c_k$ with $d_k$.
If $a_k=c_k\neq 0$ or $b_k=d_k\neq 0$ and $k$ is sufficiently large, we have $\gamma_1,\gamma_2>0$ and $\gamma_1,\gamma_2\to 0$, for $k\to\infty$.
\endthm\rm
\section{Outlook: The bifurcation problem}
The results of Theorem~\ref{thm_main} motivate to study bifurcation for the problem~\eqref{problem}; cf.~\cite{ZECui08b} where bifurcation for the non-necrotic version of our strip-shaped tumor growth model is established. In some older papers, bifurcation solutions for radially symmetric free boundary value problems have been constructed by using a power series technique; see, e.g., \cite{FR00,FHV01}. The modern method of analysis in \cite{ZECui08b} is based on an application of the following theorem of Crandall and Rabinowitz.
\thm[see \cite{CR71}] Let $X$ and $Y$ be real Banach spaces and let $G(u,\lambda)$ be a $C^q$ map ($q\geq 3$) from a neighborhood of a point $(u_0,\lambda_0)\in X\times\R$ into $Y$. We assume that
\begin{enumerate}
\item $G(u_0,\lambda_0)=G_\lambda(u_0,\lambda_0)=0$,
\item $\text{\rm Ker } G_u(u_0,\lambda_0)$ is one-dimensional and spanned by $u_0$,
\item $\text{\rm Im }G_u(u_0,\lambda_0)$ has codimension $1$ and
\item $G_{\lambda\lambda}(u_0,\lambda_0)\in \text{\rm Im }G_u(u_0,\lambda_0)$, $G_{u\lambda}(u_0,\lambda_0)u_0\notin \text{\rm Im }G_u(u_0,\lambda_0).$
\end{enumerate}
Then $(u_0,\lambda_0)$ is a bifurcation point of the equation $G(u,\lambda)=0$ in the sense that in a neighborhood of $(u_0,\lambda_0)$, the set of solutions of $G(u,\lambda)=0$ consists of two $C^{q-2}$ smooth curves $\Upsilon_{1,2}$ which intersect only at the point $(u_0,\lambda_0)$ and can be parameterized as follows:
$$
\begin{array}{lllll}
\Upsilon_1: & (u(\lambda),\lambda),    & |\lambda-\lambda_0|\text{ is small,} & u(\lambda_0)=u_0, & u'(\lambda_0)=0, \\
\Upsilon_2: & (u(\eps),\lambda(\eps)), & |\eps|\text{ is small,} & (u(0),\lambda(0))=(u_0,\lambda_0), & u'(0)=u_0.
\end{array}
$$
\endthm\rm
In \cite{ZECui08b}, the authors let the surface tension coefficient $\gamma$ (which corresponds to $\gamma_2$ in our model) play the role of the bifurcation parameter $\lambda$. Our crucial problem is that there are two surface tension coefficients $\gamma_1\neq\gamma_2$ in the necrotic variant of the multi-layer tumor growth model which suggests that the Crandall-Rabinowitz Theorem is not suitable for our purposes. Moreover, the technique presented in \cite{ZECui08b} is fairly standard to obtain bifurcation branches for related free boundary models and has already been applied in various publications \cite{FH07,EM11}; see also \cite{HHHLSZ} where bifurcation from radially symmetric solutions of a necrotic tumor growth model is discussed. Because of that, we provide some supplementary material concerning the model \eqref{problem} in this outlook and prepare a functional analytic formulation which might be suitable to derive bifurcation. It remains an open problem to establish the existence of bifurcation branches for which we probably need some deep and new ideas.

Let $h^{m+\alpha}(\S)$, $m\in\{2,3,\ldots\}$ and $\alpha\in(0,1)$, denote the little H\"older space on the circle, i.e., the closure of $C^{\infty}(\S)$ in the H\"older space $C^{m+\alpha}(\S)$. The cone of positive functions in $h^{m+\alpha}(\S)$ is denoted as $h_+^{m+\alpha}(\S)$. We define $h^{m+\alpha}(\overline\Omega)$ analogously. The small H\"older spaces are used frequently since they are Banach algebras (under pointwise multiplication) and the embedding $h^r(\S)\hookrightarrow h^s(\S)$, $r>s$, is compact. To keep our notation as simple as possible, we will label the coordinates in $\Omega$ by $x$ and $y$ in the sequel.

First, we establish that the linear second-order differential operator
$$\A(\rho_1,\rho_2)=a_{11}\partial_{x}^2+2a_{12}\partial_x\partial_y+a_{22}\partial_y^2+b\partial_y,\quad\rho_1,\rho_2\in h^{m+\alpha}_+(\S),$$
introduced in Section~\ref{sec_3} is uniformly elliptic in $\Omega$ as defined in \cite{GT77}, p.~30. Pick $\xi=(\xi_1,\xi_2)\in\R^2\backslash\{0\}$. Since the coefficients $a_{ij}=a_{ji}$, for $i,j=1,2$, are continuous functions on $\S\times(0,1)$, it is clear that $\sum_{i,j}a_{ij}\xi_i\xi_j\leq\Lambda|\xi|^2$, for some $\Lambda>0$. On the other hand, we have
$$\sum_{i,j=1}^2a_{ij}\xi_i\xi_j=\left(\xi_1-\frac{y(\rho_2'-\rho_1')+\rho_1'}{\rho_2-\rho_1}\xi_2\right)^2+\frac{1}{(\rho_2-\rho_1)^2}\xi_2^2=|\tilde\xi|^2,$$
where
$$\tilde\xi=
\left(
  \begin{array}{cc}
    1 & -\frac{y(\rho_2'-\rho_1')+\rho_1'}{\rho_2-\rho_1} \\
    0 & \frac{1}{\rho_2-\rho_1} \\
  \end{array}
\right)
\xi
=\mathbf{A}\xi.
$$
Clearly, $\mathbf A$ is invertible and there is $\lambda>0$ such that
$$\sum_{i,j=1}^2a_{ij}\xi_i\xi_j\geq\frac{1}{|\mathbf{A}^{-1}|_F^2}|\xi|^2\geq\lambda|\xi|^2;$$
here, $|\mathbf{A}^{-1}|_F$ is the Frobenius norm of $\mathbf A^{-1}$.

For any given $\rho_1,\rho_2\in h^{m+\alpha}_+(\S)$ we solve the elliptic boundary value problem
\bea
\left\{
\begin{array}{rcll}
  \A(\rho_1,\rho_2)\sigma & = & \sigma & \text{in }\Omega, \\
  \B_1(\rho_1,\rho_2)\sigma & = & 0 & \text{on }\S, \\
  \sigma|_{y=1} & = & \bar\sigma & \text{on }\S,
\end{array}
\right.
\label{3.1}
\eea
and obtain a unique solution $\sigma\in h^{m+\alpha}(\overline\Omega)$. We let $\RR$ denote the solution operator for the nutrient concentration and write $\sigma=\RR(\rho_1,\rho_2)\bar\sigma$. It follows from elliptic regularity theory \cite{GT77} that
$$\RR(\cdot,\cdot)\bar\sigma\in C^{\infty}(h_+^{m+\alpha}(\S)^2,h^{m+\alpha}(\overline\Omega)).$$
Second, we study the Neumann problem
\bea
\left\{
\begin{array}{rcll}
  \A(\rho_1,\rho_2)p & = & -\mu(\RR(\rho_1,\rho_2)\bar\sigma-\tilde\sigma) & \text{in }\Omega, \\
  \B_1(\rho_1,\rho_2)p & = & 0 & \text{on }\S, \\
  \B_2(\rho_1,\rho_2)p & = & 0 & \text{on }\S,
\end{array}
\right.
\label{3.3}
\eea
which is solvable if and only if
$$\Phi(\rho_1,\rho_2):=\int_\Omega\left(\RR(\rho_1,\rho_2)\bar\sigma-\tilde\sigma\right)(\rho_2-\rho_1)\d x\d y=0.$$
Here, $\Phi\in C^{\infty}(h_+^{m+\alpha}(\S)^2,\R)$. Using \eqref{sigmastar} and \eqref{condition1}, we compute
\bea\Phi(\rho_{1*},\rho_{2*}) \!\! & = & \!\! \delta\int_\Omega\left(\RR(\rho_{1*},\rho_{2*})\bar\sigma-\tilde\sigma\right)\d x\d y \nonumber\\
& = & \!\! \delta\bar\sigma\int_\Omega\left[\frac{\cosh(y\delta+\rho_{1*})-\tanh\rho_{1*}\sinh(y\delta+\rho_{1*})}{\cosh\rho_{2*}-\tanh\rho_{1*}\sinh\rho_{2*}}
-\frac{\tanh\delta}{\delta}\right]\dw x\d y\nonumber\\
%
%& = & \!\! \bar\sigma\left[\frac{\sinh\rho_{2*}-\sinh\rho_{1*}-\tanh\rho_{1*}(\cosh\rho_{2*}-\cosh\rho_{1*})}
%{\cosh\rho_{2*}-\tanh\rho_{1*}\sinh\rho_{2*}}-\tanh\delta\right]\nonumber\\
%
%& = & \!\! \bar\sigma\left[\frac{\tanh\rho_{2*}-\tanh\rho_{1*}}{1-\tanh\rho_{1*}\tanh\rho_{2*}}-\tanh\delta\right] \nonumber\\
%
& = & \!\! 0\nonumber\\\label{phi=0}
\eea
and, in view of \eqref{phi=0},
\bea
&& \hspace{-1cm}\Phi'(\rho_{1*},\rho_{2*})(0,1) \nonumber\\
& = & \lim_{\eps\to 0}\frac{1}{\eps}\left[\Phi(\rho_{1*},\rho_{2*}+\eps)-\Phi(\rho_{1*},\rho_{2*})\right]\nonumber\\
& = & \lim_{\eps\to 0}\frac{1}{\eps}\int_\Omega\left[\RR(\rho_{1*},\rho_{2*}+\eps)\bar\sigma-\tilde\sigma\right](\delta+\eps)\d x\d y\nonumber\\
& = & \lim_{\eps\to 0}\frac{\delta}{\eps}\int_\Omega\left[\RR(\rho_{1*},\rho_{2*}+\eps)\bar\sigma-\tilde\sigma\right]\dw x\d y\nonumber\\
& = & \lim_{\eps\to 0}\left[\frac{\delta}{\eps}\int_\Omega\left[\RR(\rho_{1*},\rho_{2*}+\eps)\bar\sigma-\RR(\rho_{1*},\rho_{2*})\bar\sigma\right]\dw x\d y
+\frac{1}{\eps}\Phi(\rho_{1*},\rho_{2*})\right]\nonumber\\
& = & \delta\int_\Omega\lim_{\eps\to 0}\frac{1}{\eps}\left[\RR(\rho_{1*},\rho_{2*}+\eps)\bar\sigma-\RR(\rho_{1*},\rho_{2*})\bar\sigma\right]\dw x \d y \nonumber\\
& = & \delta\bar\sigma\int_\Omega\left.\frac{\dw}{\dw\eps}\right|_{\eps=0}\frac{\cosh(y(\delta+\eps)+\rho_{1*})-
\tanh\rho_{1*}\sinh(y(\delta+\eps)+\rho_{1*})}{\cosh(\rho_{2*}+\eps)-\tanh\rho_{1*}\sinh(\rho_{2*}+\eps)}\d x\d y\nonumber\\
& = & \bar\sigma\left(1-\frac{1}{\delta}\tanh\delta-\tanh^2\delta\right);\nonumber\\\label{phi'neq0}\eea
the explicit calculations in the last steps are left to the reader. The expression obtained in \autoref{phi'neq0} is nonzero, since $1-\frac{1}{\delta}\tanh\delta-\tanh^2\delta=0$ implies that $\delta=\sinh\delta\cosh\delta$ which is possible only if $\delta=0$.

It follows from \eqref{phi'neq0} and the continuity of $\Phi'$ that there is a neighborhood $U$ of $(\rho_{1*},\rho_{2*})$ in $h_+^{m+\alpha}(\S)^2$ such that $\Phi'(\rho_1,\rho_2)(0,1)\neq 0$ for all $(\rho_1,\rho_2)\in U$. Thus $\Phi(\rho_1,\rho_2)=0$ defines a smooth Banach submanifold $\M$ of codimension $1$ in a small neighborhood of $(\rho_{1*},\rho_{2*})$, i.e.,
$$\M=\set{(\rho_1,\rho_2)\in h_+^{m+\alpha}(\S)^2}{\max_{i=1,2}\norm{\rho_i-\rho_{i*}}_{h_+^{m+\alpha}(\S)}<\delta,\;\Phi(\rho_1,\rho_2)=0}.$$

For any $(\rho_1,\rho_2)\in\M$ the problem~\eqref{3.3} has a solution which is unique up to a constant. Let $\T(\rho_1,\rho_2)$ be the solution operator for \eqref{3.3} which associates to the right-hand side $-f$ the solution $p$ which is zero at the origin. We then have
$$p=\mu\T(\rho_1,\rho_2)\left(\RR(\rho_{1},\rho_{2})\bar\sigma-\tilde\sigma\right)$$
and
$$\T\in C^{\infty}(h_+^{m+\alpha}(\S)^2,\L(h^{m+\alpha}(\overline\Omega),h^{m+2+\alpha}(\overline\Omega)).$$
Let $\text{tr}_i$, $i=1,2$, denote the trace operator for $\Gamma_1$ and $\Gamma_2$ respectively. We now set
$$\calS_i(\rho_1,\rho_2)=\mu\,\text{tr}_i\circ\T(\rho_1,\rho_2)\left(\RR(\rho_{1},\rho_{2})\bar\sigma-\tilde\sigma\right)$$
and recall that the curvature operator is given by
$$\kappa(\rho)=-\frac{\partial^2\rho}{\partial x^2}\left[1+\left(\frac{\partial\rho}{\partial x}\right)^2\right]^{-3/2}.$$
We have shown that the problem~\eqref{trsf} can be rewritten as
\beq\label{3.8}
\left\{
\begin{array}{rclcl}
  \calS_1(\rho_1,\rho_2) & + & c_0 + c & = & \gamma_1\kappa(\rho_1), \\
  \calS_2(\rho_1,\rho_2) & + & c_0     & = & \gamma_2\kappa(\rho_2), \\
  &&\Phi(\rho_1,\rho_2)  & = & 0,
\end{array}
\right.
\eeq
where $\calS_i\in C^{\infty}(\M,h^{m+\alpha}(\S))$ and $\kappa\in C^{\infty}(\M,h^{m-2+\alpha}(\S))$. With $\tilde\calS_1=\calS_1+c+c_0$ and $\tilde\calS_2=\calS_2+c_0$, we conclude that \eqref{3.8} is equivalent to
\beq\label{F}
F(\rho_1,\rho_2,\gamma_1,\gamma_2):=
\left(
  \begin{array}{c}
    \tilde\calS_1(\rho_1,\rho_2)-\gamma_1\kappa(\rho_1) \\
    \tilde\calS_2(\rho_1,\rho_2)-\gamma_2\kappa(\rho_2) \\
  \end{array}
\right)
=
\left(
  \begin{array}{c}
    0 \\
    0 \\
  \end{array}
\right)
\eeq
and $F\in C^{\infty}(\M\times\R_+^2,h^{m-2+\alpha}(\S)^2)$.

While the non-necrotic model in \cite{ZECui08b} can be rewritten as the zero level set of a function $X\times\R_+\to Y$, where $X$ and $Y$ are real Banach spaces suitable for the application of the Crandall-Rabinowitz Theorem, the mapping $F$ in \eqref{F} employs two positive parameters $\gamma_1,\gamma_2$ as inputs which is not compatible with the assumptions of Crandall-Rabinowitz. Thus the bifurcation problem for \eqref{F} remains a subject for further research.
\end{document}